\documentclass[11pt]{article}
\usepackage{geometry}
\usepackage{algorithm}
\usepackage{algorithmic}
\usepackage{amsthm}
\usepackage{amsmath}
\usepackage{amssymb}
\usepackage{epstopdf}
\usepackage{cite, bm, braket}
\usepackage[colorlinks,linkcolor=blue,citecolor=blue,anchorcolor=blue]{hyperref}

\title{Quantum radial basis function methods for scattered data fitting
	\thanks{This work was supported by the National Natural Science Foundation of China under Grant No.\ 11571265 and NSFC-RGC No.\ 11661161017.}}

\author{Lingxia Cui\footnotemark[2]\
	\and Hua Xiang\footnotemark[2] \footnotemark[3] \footnotemark[4]}

\date{}
\begin{document}
	\maketitle
	\renewcommand{\thefootnote}{\fnsymbol{footnote}}
	\footnotetext[2]{School of Mathematics and Statistics, Wuhan University, Wuhan 430072, China.}
	\footnotetext[3]{Hubei Key Laboratory of Computational Science, Wuhan University, Wuhan 430072, China.}
	\footnotetext[4]{Corresponding author. E-mail address: {\tt hxiang@whu.edu.cn}.}
\begin{abstract}
Scattered data fitting is a frequently encountered problem for reconstructing an unknown function from given scattered data. Radial basis function (RBF) methods have proven to be highly useful to deal with this problem. We describe two quantum algorithms to efficiently fit scattered data based on globally and compactly supported RBFs respectively. For the globally supported RBF method, the core of the quantum algorithm relies on using coherent states to calculate the radial functions and a nonsparse matrix exponentiation technique for efficiently performing a matrix inversion. A quadratic speedup is achieved in the number of data over the classical algorithms. For the compactly supported RBF method, we mainly use the HHL algorithm as a subroutine to design an efficient quantum procedure that runs in time logarithmic in the number of data, achieving an exponential improvement over the classical methods.

\end{abstract}

\section{Introduction}
Data fitting is a practical problem that has many applications in science and engineering. In many cases the data sites are scattered, bearing no regular structure at all or existing in a space of very high dimensions. Many traditional numerical methods are usually limited to very special situations and cannot efficiently handle such problems except the RBF methods, also known as meshfree methods. The RBF methods are often better suited to cope with the problems involving high space dimensions or change in the complex geometry of the domain of interest. Its applications can be found in many different areas, such as terrain modeling \cite{hardy1971multiquadric}, surface reconstruction \cite{morse2001interpolating}, fluid-structure interaction \cite{beckert2001multivariate}, the numerical solution of partial differential equations \cite{franke1998solving, kansa2000circumventing, wu2004dynamically}, kernel learning \cite{scholkopf2002learning, schaback2006kernel}, and parameter estimation, to just name a few. In this paper, we will focus on the multivariate scattered data interpolation problem, which is one of the most important and fundamental problems in data fitting and approximation theory in general. For RBF methods, a linear combination of certain RBFs which may be globally or compactly supported, is often used to interpolate given data. This further leads to a dense or sparse system of linear equations. Hence solving the interpolation problem amounts to obtaining the solution of the system of linear equations. However, if the number of data is large, the dimension of the linear system may be high and the fitting can become very costly.

With the rise of quantum information theory, finding a quantum algorithm that efficiently solves this problem is possible. By exploiting quantum mechanical effects, large-scale research aimed at finding quantum algorithms which are more efficient than the best classical counterparts has been set up \cite{harrow2009quantum, wiebe2012quantum, lloyd2014quantum, rebentrost2014quantum, wang2019quantum, shao2020data}. Recently, the breakthrough work of Harrow, Hassidim and Lloyd \cite{harrow2009quantum} introduced a quantum linear system algorithm, called HHL algorithm, that computes the quantum state corresponding to the solution of a linear system for a sparse and well-conditioned matrix, in time exponentially faster than the best classical algorithm. Subsequent improvements of the HHL algorithm for solving systems of linear equations can be found in \cite{ambainis2010variable, childs2017quantum, wossnig2018quantum} and so on. The quantum linear solver lies at the heart of many quantum algorithms and has inspired several works in the area of quantum data fitting, such as \cite{wiebe2012quantum, liu2017fast, wang2017quantum, shao2020quantum}. In \cite{wiebe2012quantum}, a quantum data-fitting algorithm was proposed for the least squares method. In \cite{liu2017fast}, the authors gave efficient quantum algorithms for least squares regression and statistic leverage scores. A quantum algorithm which can output the optimal parameters in the classical form was developed in \cite{wang2017quantum} for linear regression by the least squares approach. The authors in \cite{shao2020quantum} considered quantum regularized least squares solver with parameter estimate. For the total least squares method, a quantum data-fitting algorithm based on quantum resonant transitions was presented in \cite{wang2019quantum}. Different from these listed methods above, we consider using RBF methods to solve the scattered data-fitting problem.

In this paper, we show that a quantum scattered data-fitting algorithm can be implemented with an improved performance over the classical counterpart. We describe two quantum methods using globally and locally supported RBFs respectively. For the globally supported RBF method, we focus on using the Gaussian RBF, which is one of the most important global RBFs and plays a central role in many applications, to design a quantum algorithm for a well-posed scattered data interpolation problem. The key to this quantum method is the use of coherent states \cite{sanders2012review, chatterjee2016generalized} and a matrix exponentiation technique \cite{lloyd2014quantum}. With the help of the coherent states, we can efficiently calculate all pairs of the data sites to construct the coefficient matrix of the linear system in time $\widetilde O(md)$ with $m$ the number of data and $d$ the dimension of the input vectors. The $\widetilde{O}$ notation suppresses more slowly growing factors. Since the coefficient matrix is a nonsparse density matrix, we then employ the technique for the exponentiation of nonsparse matrices developed in \cite{lloyd2014quantum} to solve the linear system and obtain a quantum state of the solution in time which scales linearly in $m$. A quadratic speedup is achieved in the number of data over the classical methods. For the compactly supported RBF method, we mainly use the HHL algorithm as a subroutine to deal with the problem, due to the sparsity of the coefficient matrix. Similar to the globally supported RBF method, we first design an oracle to calculate all the nonzero entries of the coefficient matrix, then perform a matrix inversion algorithm based on the Hamiltonian simulation subroutine in time that scales logarithmically in $m$ and $d$, exponentially faster than the classical algorithms. For the evaluation on a new given data, we apply a swap test to estimate the value of the new data.

The outline of the paper is as follows. In Section 2, some preliminaries are introduced. In Section 3, we state the problem of scattered data interpolation using RBF methods. The main results are presented in Section 4, where we give two quantum RBF methods for the scattered data interpolation problem based on globally and compactly supported RBFs respectively. Finally, Section 5 provides a brief conclusion and a discussion of future work.

\section{Preliminaries}
For any matrix $A\in \mathbb{R}^{m\times m}$, its Frobenius norm is defined by $||A||_F = \sqrt{\sum_{i,j=1}^mA_{ij}^2}$ with $A_{ij}$ denoting the $(i,j)$th element of $A$. Its 2-norm is defined by $||A||_2 = \sigma_{\max}$, and its 2-norm condition number is given by
$$
\kappa(A)=||A||_2||A^{-1}||_2=\frac{\sigma_{\max}}{\sigma_{\min}}
$$
with the convention that $\kappa=\infty$ for the singular matrix $A$. Here $\sigma_{\max}$ and $\sigma_{\min}$ are the largest and smallest singular values of $A$ respectively. If $A$ is a symmetric positive definite matrix, then the condition number of $A$ can also be computed as the ratio of the largest and smallest eigenvalues, i.e.,
$$
\kappa=\frac{\lambda_{\max}}{\lambda_{\min}}.
$$
For a vector $\bm{x}=[x_1,x_2,\dots,x_d]^T\in\mathbb{R}^d$, we denote the Euclidean norm $||\bm{x}|| = (\sum_{i=1}^d {x_i}^2)^{1/2}$.  Besides, we denote $\bm{x}^{(i)}$ as the $i$th vector and ${x}_j$ as the $j$th component of vector $\bm{x}$. In the quantum setting, we use the standard bra-ket notation to denote quantum states. The vector state $\ket{\bm{x}}$ for $\bm{x}\in\mathbb{R}^d$ is defined as ${1}/{||\bm{x}||}\sum_{i=1}^d x_i \ket{i}$.

For the purpose of error analysis, two lemmas about matrix analysis are given as follows.
\newtheorem{lemma}{Lemma}
\begin{lemma}\label{lemma.1}
\cite{golub2014matrix}
If $A, E \in \mathbb{R}^{m\times m}$, $A$ is nonsingular and $r\equiv ||A^{-1}E||_p<1$, then $A+E$ is nonsingular and
$$
||(A+E)^{-1}-A^{-1}||_p \leq \frac{||E||_p||A^{-1}||_p^2}{1-r}
$$
with $||\cdot||_p$ denoting the $p-$norm ($p\geq1$).
\end{lemma}

\begin{lemma}\label{lemma.2}
\cite{golub2014matrix}
If $A$ and $A+E$ are $m$-by-$m$ symmetric matrices, then
$$
|\lambda_k(A+E)-\lambda_k(A)| \leq ||E||_2
$$
for $k=1,2,\dots,m$, where the notations $\lambda_k(A+E)$ and $\lambda_k(A)$ represent the $k$th largest eigenvalues of $A+E$ and $A$ respectively.
\end{lemma}

Furthermore, we introduce some basic concepts and properties about radial functions for completeness.
\newtheorem{definition}{Definition}
\begin{definition}
\cite{fasshauer2007meshfree}
A function $\Phi: \mathbb{R}^d \rightarrow \mathbb{R}$ is called radial provided that there exists a univariate function $\phi:[0,\infty) \rightarrow \mathbb{R}$ such that
$$
\Phi(\bm x)=\phi(r),
$$
where $r=||\bm x||$.
\end{definition}

\begin{definition}
\cite{fasshauer2007meshfree}
A complex-valued continuous function $\Phi: \mathbb{R}^d \rightarrow \mathbb{C}$ is called positive definite if
$$
\sum\limits_{i=1}^{m}\sum\limits_{j=1}^{m}c_i\bar{c}_j\Phi(\bm x^{(i)}-\bm x^{(j)})>0
$$
for any $m$ pairwise different points $\{{\bm x^{(1)},\bm x^{(2)},\dots,\bm x^{(m)}}\} \in \mathbb{R}^d$, and $\bm c=[c_1,c_2,\dots,c_m] \in \mathbb{C}^m \backslash \{\bm 0\}$.
\end{definition}

\begin{definition}
\cite{fasshauer2007meshfree}
We call a univariate function $\phi:[0,\infty)\rightarrow \mathbb{R}$ positive definite on $\mathbb{R}^d$ if the corresponding multivariate function $\Phi(\bm x)\equiv\phi(||\bm x||), \bm x\in \mathbb{R}^d$, is positive definite.
\end{definition}

\begin{lemma}\label{lemma.3}
\cite{fasshauer2007meshfree}
If $\Phi$ is positive definite on $\mathbb{R}^d$, then\\
{\rm{(i)}} $\Phi(\bm 0)\geq0$.\\
{\rm{(ii)}} $\Phi(-\bm x)=\overline{\Phi(\bm x)}$.\\
{\rm{(iii)}} $\Phi$ is bounded. In fact, $|\Phi(\bm x)|\leq \Phi(\bm 0)$.
\end{lemma}

A radial function is radially symmetric around its center $\bm x^{(c)}\in\mathbb{R}^d$. For any fixed center $\bm{x}^{(c)}$, it can be defined as $\Phi(\bm{x}-\bm{x}^{(c)}) = \phi(r)$, where $r = ||\bm{x} - \bm{x}^{(c)}||$ is the Euclidean distance. To use $\phi(r)$ as a basis function in a RBF method,  $\bm{x}$ is taken to be the input variable and the center $\bm{x}^{(c)}$ is set to a constant point. Note that the univariate function $\phi$ is independent from the dimension $d$ of the input vectors. Consequently, methods that use radial functions as their basis functions are, in principle, easily adapted to solve problems in high dimensions.

\section{Scattered data interpolation}
Scattered data interpolation can be achieved using RBFs centered at the constraints. We are now ready for a precise formulation of the scattered data interpolation problem.

\newtheorem{problem}{Problem}
\begin{problem}\label{problem.1}
Given a set of $m$ distinct data sites $\{\bm{x}^{(1)}, \bm{x}^{(2)}, \dots, \bm{x}^{(m)}\}$ in $\mathbb{R}^d$ and a corresponding set of $m$ values $\{y^{(1)}, y^{(2)},\dots, y^{(m)}\}$ in $\mathbb{R}$ sampled from some unknown function $f^{(exact)}$, i.e., $y^{(j)} = f^{(exact)}(\bm{x}^{(j)})$, $j=1,2,\dots,m$, find a function $f$ such that $f(\bm{x}^{(j)}) = y^{(j)}$ for $j=1,2,\dots,m$.
\end{problem}

A convenient and common approach for solving the scattered data problem is the RBF methods. We can choose a radial function $\phi$ and a set of centers $\{\bm{x}^{(c_1)},\bm{x}^{(c_2)}, \dots, \bm{x}^{(c_{k})}\}$ for some $k\in\mathbb{N}$, to obtain a basis $\{\phi(||\cdot-\bm{x}^{(c_1)}||), \phi(||\cdot-\bm{x}^{(c_2)}||), \dots, \phi(||\cdot-\bm{x}^{(c_k)}||)\}$. This basis can then be used to construct an approximation $f$ of the function $f^{(exact)}$.

Without loss of generality, we choose the centers to coincide with the data sites, i.e., $\bm{x}^{(c_j)} = \bm{x}^{(j)}$ for each $j = 1,2,\dots,m$. In this sense, the basis is $\{\phi(||\cdot-\bm{x}^{(1)}||), \phi(||\cdot-\bm{x}^{(2)}||),\dots, \phi(||\cdot-\bm{x}^{(m)}||)\}$.
The approximation $f$ can then be constructed from a linear combination of these $m$ RBFs, so that
$$
f(\bm{x}) = \sum\limits_{j=1}^{m}c_j \phi(||\bm{x} - \bm x^{(j)}||)
$$
with constant coefficients $c_j$. The coefficients $c_j$ are determined by ensuring that the approximation function will exactly match the given data at the data sites, i.e., $f(\bm{x}^{(j)}) = y^{(j)}$ for $j=1,2,\dots,m$, which produces a system of linear equations of the form
\begin{equation}\label{eq.1}
	\left[
	\begin{array}{cccc}
		\phi(||\bm{x}^{(1)} - \bm{x}^{(1)}||) & \phi(||\bm{x}^{(1)} - \bm{x}^{(2)}||) & \cdots & \phi(||\bm{x}^{(1)} - \bm{x}^{(m)}||)\\
		\phi(||\bm{x}^{(2)} - \bm{x}^{(1)}||) & \phi(||\bm{x}^{(2)} - \bm{x}^{(2)}||) & \cdots & \phi(||\bm{x}^{(2)} - \bm{x}^{(m)}||)\\
		\vdots & \vdots & \ddots & \vdots \\
		\phi(||\bm{x}^{(m)} - \bm{x}^{(1)}||) & \phi(||\bm{x}^{(m)} - \bm{x}^{(2)}||) & \cdots & \phi(||\bm{x}^{(m)} - \bm{x}^{(m)}||)
	\end{array}
	\right]
	\left[
	\begin{array}{c}
	c_1\\ c_2\\ \vdots\\ c_m
	\end{array}
	\right]
	=
	\left[
	\begin{array}{c}
		y^{(1)}\\ y^{(2)}\\ \vdots\\ y^{(m)}
	\end{array}
	\right].
\end{equation}
Or equivalently, we find the coefficients $c_j$ by solving the linear system
\begin{equation}\label{eq.2}
	\mathbf{A}\bm c^{(exact)} = \bm y,
\end{equation}
where
\begin{equation}
\mathbf{A} = \frac{1}{m}
\left[
\begin{array}{cccc}
	\phi(||\bm{x}^{(1)} - \bm{x}^{(1)}||) & \phi(||\bm{x}^{(1)} - \bm{x}^{(2)}||) & \cdots & \phi(||\bm{x}^{(1)} - \bm{x}^{(m)}||)\\
	\phi(||\bm{x}^{(2)} - \bm{x}^{(1)}||) & \phi(||\bm{x}^{(2)} - \bm{x}^{(2)}||) & \cdots & \phi(||\bm{x}^{(2)} - \bm{x}^{(m)}||)\\
	\vdots & \vdots & \ddots & \vdots \\
	\phi(||\bm{x}^{(m)} - \bm{x}^{(1)}||) & \phi(||\bm{x}^{(m)} - \bm{x}^{(2)}||) & \cdots & \phi(||\bm{x}^{(m)} - \bm{x}^{(m)}||)
\end{array}
\right],
\end{equation}
\begin{equation*}
\bm c^{(exact)} = \left[
\begin{array}{c}
c_1 \\ c_2 \\ \vdots \\ c_m	
\end{array}
\right], \quad
\bm y = \frac{1}{m}\left[
\begin{array}{c}
	y^{(1)} \\ y^{(2)} \\ \vdots \\ y^{(m)}	
\end{array}
\right].
\end{equation*}
For the convenience of the following analysis, we refer both the coefficient matrices in (\ref{eq.1}) and (\ref{eq.2}) as the interpolation matrices. Obviously, they are both symmetric matrices.

It is clear that the solution of the scattered data interpolation problem hinges entirely on the solution of the system of linear equations. We know from linear algebra that the system will have a unique solution whenever the interpolation matrix is nonsingular. With proper choice of the basis function such as the Gaussian RBF or inverse multiquadric RBF, the interpolation matrix can be positive definite and thus nonsingular. The coefficients can then be found and the approximation function $f$ can be constructed.

\textbf{Remark 1.} \ To make sure the uniqueness of the solution of the linear system, we select positive definite radial functions as basis functions \cite{fasshauer2007meshfree}. Besides, we assume that the problem is well-posed, that is, the interpolation matrix is well-conditioned for the numerical stability.

\section{$\mathbf{Quantum \; RBF \; methods}$}
We notice in the previous section that the solution of the scattered data interpolation problem using RBF methods boils down to solving the linear system (\ref{eq.1}) or (\ref{eq.2}). In this section, we present two quantum algorithms to deal with the problem based on globally and compactly supported RBFs respectively. We mainly adapt a strategy developed by Lloyd et al. in \cite{lloyd2014quantum} to compute the matrix inverse for the globally supported RBF method. As for the compactly supported RBF method, HHL algorithm is employed for solving the linear system. In the following, we design these algorithms and describe their computational advantage over the classical counterparts.

Before giving the main results, we assume that oracle for data sites that returns quantum vectors $\ket{\bm x^{(i)}}$, the norms $||\bm x^{(i)}||$ for $i=1,2,\dots,m$,  $\ket{\bm y}$ and $||\bm y||$ is given, such as by quantum random access memory (QRAM) \cite{giovannetti2008quantum}, which uses only $O\left(\log(md)\right)$ operations to access them.

\subsection{$\mathbf{Quantum \; globally \; supported \; RBF \; method}$}
Global interpolation methods based on RBFs can be easily implemented for fitting functions sampled at scattered data points. They appear in many different areas, such as geology, biology, engineering and statistics. Table \ref{table.1} lists some common globally supported RBFs, in which $\eta$ is a scaling parameter and the variable $r$ stands for $||\bm x-\bm x^{(c)}||$ as before. With the exception of the multiquadric, all of the functions listed are positive definite \cite{fasshauer2007meshfree}.
\begin{table}[h]
	\centering
	\caption{\label{table.1}Some common global radial basis functions.}
	\begin{tabular}{ll}\hline\hline
		Name & Definition\\ \hline
		Gaussian & $\phi(r)=e^{-(\eta r)^2}$\\
		Multiquadric & $\phi(r)=\sqrt{1+(\eta r)^2}$\\
		Inverse Multiquadric & $\phi(r)=\frac{1}{\sqrt{1+(\eta r)^2}}$\\
		$C^0$ Mat\'ern & $\phi(r)=e^{-\eta r}$\\
		$C^2$ Mat\'ern & $\phi(r)=e^{-\eta r}(1+\eta r)$\\
		$C^4$ Mat\'ern & $\phi(r)=e^{-\eta r}(3+3\eta r+(\eta r)^2)$\\
		\hline\hline		
	\end{tabular}
\end{table}

In a global RBF interpolation, Gaussian
$$
\phi(r)=e^{-r^2/{2\sigma^2}}
$$
is a good candidate for a positive definite function, where $\sigma=\sqrt{1/(2\eta^2)}$ is defined by the scaling parameter for convenience. It is one of the most important globally supported RBFs and plays a central role in different areas, especially in statistics. Under some mild conditions, Gaussian performs rather well in practice. In this subsection, we concentrate on using Gaussian RBF to design a quantum algorithm for the scattered data interpolation problem.

With Gaussian RBF, the approximation function $f$ is of the form
\begin{equation}\label{eq.4}
	f(\bm{x}) = \sum\limits_{j=1}^{m}c_j e^{-||\bm{x} - \bm x^{(j)}||^2/2\sigma^2}.
\end{equation}
The coefficients $c_j$ are determined by solving the linear system (\ref{eq.2}) with
\begin{equation}\label{eq.5}
\mathbf{A} = \frac{1}{m}
\left[
\begin{array}{cccc}
e^{-||\bm{x}^{(1)} - \bm{x}^{(1)}||^2/{2\sigma^2}} & e^{-||\bm{x}^{(1)} - \bm{x}^{(2)}||^2/{2\sigma^2}} & \cdots & e^{-||\bm{x}^{(1)} - \bm{x}^{(m)}||^2/{2\sigma^2}}\\
e^{-||\bm{x}^{(2)} - \bm{x}^{(1)}||^2/{2\sigma^2}} & e^{-||\bm{x}^{(2)} - \bm{x}^{(2)}||^2/{2\sigma^2}} & \cdots & e^{-||\bm{x}^{(2)} - \bm{x}^{(m)}||^2/{2\sigma^2}}\\
\vdots & \vdots & \ddots & \vdots \\
e^{-||\bm{x}^{(m)} - \bm{x}^{(1)}||^2/{2\sigma^2}} & e^{-||\bm{x}^{(m)} - \bm{x}^{(2)}||^2/{2\sigma^2}} & \cdots & e^{-||\bm{x}^{(m)} - \bm{x}^{(m)}||^2/{2\sigma^2}}
\end{array}
\right]
\end{equation}
representing the interpolation matrix, which is a symmetric positive definite matrix of dimension $m\times m$ and thus invertible. Denote $\Lambda_{\max}$ and $\Lambda_{\min}$ as its largest and smallest eigenvalues. Its condition number is
$$
\kappa=||\mathbf{A}||_2||\mathbf{A}^{-1}||_2=\Lambda_{\max}/\Lambda_{\min}.
$$
By Gershgorin's theorem and Lemma \ref{lemma.3}, it is easy to derive that
\begin{equation}\label{eq.6}
||\mathbf{A}||_2=\Lambda_{\max}\leq 1.
\end{equation}
Classically, for solving the linear system, we need to calculate each element of the interpolation matrix. There are $m(m-1)/2$ exponentials to evaluate in (\ref{eq.5}) and each element takes $O(d)$ time to calculate. Therefore, $O(m^2d)$ time is required to calculate this interpolation matrix. Furthermore, taking into account the cost for solving the linear system, it will take even more time for the problem. This can be improved by a quantum approach given as follows.

For a quantum algorithm dealing with the problem, we use a quantum linear solver to obtain a quantum state corresponding to the solution of the linear system and apply a swap test to evaluate on a new given data. Our method involves four parts: interpolation matrix preparation, interpolation matrix exponentiation, quantum algorithm implementation,  error analysis and runtime estimate.

\subsubsection{$\mathbf{Interpolation \; matrix \; preparation}$}
We note that a significant portion of the computational cost involved in the classical algorithms is the cost required to construct the interpolation matrix (\ref{eq.5}), in which all pairs of data sites need to be evaluated. In the quantum setting, we use coherent states to calculate these elements and then prepare the interpolation matrix. We now demonstrate how this is possible.

Coherent states play an important role in quantum optics and mathematical physics. They are defined in the Fock states $\{\ket{0},\ket{1},\dots\}$, which is a basis of infinite-dimensional Hilbert space. Denote $a$, $a^{\dagger}$ as the annihilation and creation operators of the harmonic oscillator respectively, that is,
$$
a \ket{k} = \sqrt{k}\ket{k-1}, \quad a^{\dagger}\ket{k}=\sqrt{k+1}\ket{k+1}
$$
with $a\ket{0}=0$. Using the creation operator, the excited state $\ket{k}$ for any $k\geq1$ can be calculated. This gives
\begin{equation}\label{eq.7}
\ket{k} = \frac{(a^\dagger)^k}{\sqrt{k!}}\ket{0}.
\end{equation}
For a real number $r\in\mathbb{R}$, its coherent state is defined by
$$
\ket{\psi^{(exact)}_r}=e^{-r^2/2\sigma^2}\sum\limits_{k=0}^{\infty}\frac{(r/\sigma)^k}{\sqrt{k!}}\ket{k}.
$$
It is a unit eigenvector of the annihilation operator $a$ corresponding to the eigenvalue $r/\sigma$, that is,
$$
a\ket{\psi^{(exact)}_r}=(r/\sigma)\ket{\psi^{(exact)}_r}.
$$
Using Eq.(\ref{eq.7}) for $\ket{k}$ yields that
$$
\ket{\psi_r^{(exact)}}=e^{-r^2/2\sigma^2}e^{ra^{\dagger}/\sigma}\ket{0}=e^{r(a^\dagger-a)/\sigma}\ket{0}.
$$
Hence, the coherent state $\ket{\psi^{(exact)}_r}$ can be generated by the unitary displacement operator $e^{r(a^\dagger-a)/\sigma}$ of dimension infinity operating on the ground state $\ket{0}$ of the harmonic oscillator.

To prepare $\ket{\psi_r^{(exact)}}$ in a finite quantum circuit, we consider its Taylor approximation
\begin{equation}\label{eq.8}
\ket{{\psi}_r} \propto \sum\limits_{k=0}^{N-1}\frac{(r/\sigma)^k}{\sqrt{k!}}\ket{k}.
\end{equation}
It is easy to obtain the inequality
\begin{equation}\label{eq.9}
\sum\limits_{k=N}^{\infty}\frac{(r/\sigma)^{2k}}{k!} \leq \frac{(r/\sigma)^{2N}}{N!}e^{(r/\sigma)^2}.
\end{equation}
For notational simplicity we define $\mathcal{B}=e^{(r/\sigma)^2}$ and $B=\sum_{k=0}^{N-1}\frac{(r/\sigma)^{2k}}{k!}$ as the normalization factors of $\ket{\psi^{(exact)}_r}$ and $\ket{\psi_r}$ respectively. Obviously, $B$ is the sum of the first $N$ terms of $\mathcal{B}$'s Taylor expansion, and $\Delta B=\mathcal{B}-B$ is the truncation error bounded by Eq.(\ref{eq.9}). Then we have
\begin{equation}\label{eq.10}
\begin{aligned}
\big|\big|\ket{\psi_r^{(exact)}}-\ket{\psi_r}\big|\big|^2	&=\Big|\Big|\mathcal{B}^{-1/2}\sum\limits_{k=0}^{\infty}\frac{(r/\sigma)^k}{\sqrt{k!}}\ket{k}-B^{-1/2}\sum\limits_{k=0}^{N-1}\frac{(r/\sigma)^k}{\sqrt{k!}}\ket{k}\Big|\Big|^2\\
	&=\Big|\Big|\left(\mathcal{B}^{-1/2}-B^{-1/2}\right)\sum\limits_{k=0}^{N-1}\frac{(r/\sigma)^k}{\sqrt{k!}}\ket{k}+\mathcal{B}^{-1/2}\sum\limits_{k=N}^{\infty}\frac{(r/\sigma)^k}{\sqrt{k!}}\ket{k}\Big|\Big|^2\\
	&=\left(\mathcal{B}^{-1/2}-B^{-1/2}\right)^2\sum\limits_{k=0}^{N-1}\frac{(r/\sigma)^{2k}}{k!} + \mathcal{B}^{-1}\sum\limits_{k=N}^{\infty}\frac{(r/\sigma)^{2k}}{k!} \\
	& \leq \left(\sqrt{\frac{B}{\mathcal{B}}}-1\right)^2+\frac{(r/\sigma)^{2N}}{N!} \\
	& \leq \left(\frac{\Delta B}{\mathcal{B}}\right)^2+\frac{(r/\sigma)^{2N}}{N!} \\
	& \leq \frac{2(r/\sigma)^{2N}}{N!}.	
\end{aligned}
\end{equation}
In Eq.(\ref{eq.10}), the first inequality is obtained from Eq.(\ref{eq.9}), the second uses the variant that $1-\sqrt{B/\mathcal{B}}=1-\sqrt{1-\Delta B/\mathcal{B}}$ and the fact that $1-\sqrt{1-z}=z/(1+\sqrt{1-z})\leq z$ for any $0\leq z\leq 1$, and the last uses that $(\Delta B/\mathcal{B})^2 \leq \Delta B/\mathcal{B}$ due to $0\leq \Delta B/\mathcal{B} \leq 1$. Set the upper bound of Eq.(\ref{eq.10}) as $\delta^2$. By Stirling's
approximation $N!\approx \sqrt{2\pi N}(N/e)^{N}$, we have
$$
\frac{2(r/\sigma)^{2N}}{\sqrt{2\pi N}(N/e)^{N}}\leq \delta^2.
$$
Taking the logarithm on both sides yields
$$
2\log\frac{1}{\delta}-\log\sqrt{2\pi}+\log 2
\leq \left(N+\frac{1}{2}\right)\log N-N\left(\log e+2\log\frac{r}{\sigma}\right).
$$
So we can choose $N=O\left(\log(1/\delta)\right)$ in Eq.(\ref{eq.8}) such that $\big|\big|\ket{\psi^{(exact)}_r}-\ket{\psi_r}\big|\big|\leq \delta$. That is, we can truncate the skew-Hermitian matrix $a^\dagger-a$ to order $O(\log 1/\delta)$, which is a two-sparse matrix, to efficiently perform the Hamiltonian simulation. Thus, there is an efficient quantum algorithm to prepare the coherent state $\ket{\psi^{(exact)}_r}$ in time $O\left(\log(1/\delta)\right)$ within error $\delta$.

The coherent state for a real vector $\bm x=(x_1,x_2,\dots,x_d)$ is defined by
\begin{equation}
\ket{\psi^{(exact)}_{\bm x}}=\ket{\psi^{(exact)}_{x_1}}\otimes\ket{\psi^{(exact)}_{x_2}}\otimes\cdots\otimes\ket{\psi^{(exact)}_{x_d}}.
\end{equation}
For any two real vectors $\bm x, \bm y\in\mathbb{R}^d$, we have
$$
\braket{\psi^{(exact)}_{\bm x}|\psi^{(exact)}_{\bm y}}=e^{-||\bm x-\bm y||^2/2\sigma^2},
$$
which is just of the form of Gaussian radial functions. Let
$$
\ket{\psi_{\bm x}}=\ket{\psi_{ x_1}}\otimes\ket{\psi_{ x_2}}\otimes\cdots\otimes\ket{\psi_{x_d}}
$$
be an approximation of $\ket{\psi^{(exact)}_{\bm x}}$, in which $\ket{\psi_{ x_i}}$ for $i=1,2,\dots,d$ are defined by Eq.(\ref{eq.8}). The error for preparing $\ket{\psi^{(exact)}_{\bm x}}$ is bounded by $d\delta$, that is,
\begin{equation}\label{eq.12}
\begin{aligned}
\big|\big|\ket{\psi^{(exact)}_{\bm x}}-\ket{\psi_{\bm x}}\big|\big|
& = \big|\big|\ket{\psi^{(exact)}_{x_1}}\otimes\ket{\psi^{(exact)}_{x_2}}\otimes\cdots\otimes\ket{\psi^{(exact)}_{x_d}}-\ket{\psi_{ x_1}}\otimes\ket{\psi_{ x_2}}\otimes\cdots\otimes\ket{\psi_{x_d}}\big|\big|\\
& \leq  \big|\big|\ket{\psi^{(exact)}_{x_1}}\otimes\ket{\psi^{(exact)}_{x_2}}\otimes\cdots\otimes\ket{\psi^{(exact)}_{x_d}}-\ket{\psi_{ x_1}}\otimes\ket{\psi^{(exact)}_{x_2}}\otimes\cdots\otimes\ket{\psi^{(exact)}_{x_d}}\big|\big|\\
& \quad +\big|\big|\ket{\psi_{ x_1}}\otimes\ket{\psi^{(exact)}_{x_2}}\otimes\cdots\otimes\ket{\psi^{(exact)}_{x_d}}-\ket{\psi_{ x_1}}\otimes\ket{\psi_{x_2}}\otimes\cdots\otimes\ket{\psi_{x_d}}\big|\big|\\
& \leq \delta + \big|\big|\ket{\psi_{ x_1}}\otimes\ket{\psi^{(exact)}_{x_2}}\otimes\cdots\otimes\ket{\psi^{(exact)}_{x_d}}-\ket{\psi_{ x_1}}\otimes\ket{\psi_{ x_2}}\otimes\cdots\otimes\ket{\psi_{x_d}}\big|\big|\\
& \leq \cdots \leq d\delta,
\end{aligned}
\end{equation}
where the first inequality applies the triangle inequality principle about vector norm, the second is based on the basic property of the Kronecker product and the last is derived by induction. It follows from Eq.(\ref{eq.12}) that
\begin{equation}\label{eq.13}
\begin{aligned}
\big|\braket{\psi^{(exact)}_{\bm x}|\psi^{(exact)}_{\bm y}}-\braket{\psi_{\bm x}|\psi_{\bm y}}\big|
& = \big|\braket{\psi^{(exact)}_{\bm x}|\psi^{(exact)}_{\bm y}}-\braket{\psi_{\bm x}|\psi^{(exact)}_{\bm y}}+\braket{\psi_{\bm x}|\psi^{(exact)}_{\bm y}}-\braket{\psi_{\bm x}|\psi_{\bm y}}\big|\\
& \leq \big|\braket{\psi^{(exact)}_{\bm x}|\psi^{(exact)}_{\bm y}}-\braket{\psi_{\bm x}|\psi^{(exact)}_{\bm y}}\big|+\big|\braket{\psi_{\bm x}|\psi^{(exact)}_{\bm y}}-\braket{\psi_{\bm x}|\psi_{\bm y}}\big|\\
& = \Big|\left(\bra{\psi^{(exact)}_{\bm x}}-\bra{\psi_{\bm x}}\right)\ket{\psi^{(exact)}_{\bm y}}\Big|+\Big|\bra{\psi_{\bm x}}\left(\ket{\psi^{(exact)}_{\bm y}}-\ket{\psi_{\bm y}}\right)\Big|\\
& \leq \big|\big|\bra{\psi^{(exact)}_{\bm x}}-\bra{\psi_{\bm x}}\big|\big|\cdot\big|\big|\ket{\psi^{(exact)}_{\bm y}}\big|\big|+\big|\big|\bra{\psi_{\bm x}}\big|\big|\cdot\big|\big|\ket{\psi^{(exact)}_{\bm y}}-\ket{\psi_{\bm y}}\big|\big|\\
& \leq 2d\delta
\end{aligned}
\end{equation}
for any two real vectors $\bm x, \bm y\in\mathbb{R}^d$.

Consider the superposition of coherent states of the data sites as follows.
$$
\ket{\Psi^{(exact)}}=\frac{1}{\sqrt{m}}\sum\limits_{j=1}^{m}\ket{j}\ket{\psi^{(exact)}_{\bm x^{(j)}}},
$$
where $\bm x^{(j)}=\left(x_1^{(j)},x_2^{(j)},\dots,x_d^{(j)}\right)$ is the vector of the $j$th data site and $\ket{\psi^{(exact)}_{\bm x^{(j)}}}=\ket{\psi^{(exact)}_{x_1^{(j)}}}\otimes\ket{\psi^{(exact)}_{x_2^{(j)}}}\otimes\cdots\otimes\ket{\psi^{(exact)}_{x_d^{(j)}}}$. It is easy to obtain the interpolation matrix (\ref{eq.5}) by using the partial trace, that is,
\begin{equation}\label{eq.14}
\begin{aligned}
{\rm tr}_2 \ket{\Psi^{(exact)}}\bra{\Psi^{(exact)}}
&=\frac{1}{m}\sum_{i,j=1}^{m}\braket{\psi^{(exact)}_{\bm x^{(i)}}|\psi^{(exact)}_{\bm x^{(j)}}}\ket{i}\bra{j}\\
&=\frac{1}{m}\sum_{i,j=1}^{m}e^{-||\bm x^{(i)}-\bm x^{(j)}||^2/{2\sigma^2}}\ket{i}\bra{j}\\
&= \mathbf{A}.
\end{aligned}
\end{equation}
In terms of a finite quantum circuit, let
$$
\ket{\Psi}=\frac{1}{\sqrt{m}}\sum\limits_{j=1}^{m}\ket{j}\ket{\psi_{\bm x^{(j)}}}
$$
denote a perturbed version of $\ket{\Psi^{(exact)}}$ with $\ket{\psi_{\bm x^{(j)}}}=\ket{\psi_{ x_1^{(j)}}}\otimes\ket{\psi_{ x_2^{(j)}}}\otimes\cdots\otimes\ket{\psi_{x_d^{(j)}}}$, the approximation of $\ket{\psi^{(exact)}_{\bm x^{(j)}}}$. It can be generated by first preparing $\ket{\psi_{\bm x^{(j)}}}$ for each $j=1,2,\dots,m$, then generating their superposition. Taking partial trace over the second register of $\ket{\Psi}\bra{\Psi}$ gives rise to the density operator of a nearby interpolation matrix
\begin{equation}\label{eq.15}
\begin{aligned}
A
& \equiv {\rm tr}_2 \ket{\Psi}\bra{\Psi}\\
& = \frac{1}{m}\sum\limits_{i,j=1}^{m}\braket{\psi_{\bm x^{(i)}}|\psi_{\bm x^{(j)}}}\ket{i}\bra{j} \\
& = \frac{1}{m}\sum\limits_{i,j=1}^{m}\big{[}\braket{\psi^{(exact)}_{\bm x^{(i)}}|\psi^{(exact)}_{\bm x^{(j)}}}+\Delta _{ij}\big{]}\ket{i}\bra{j} \\
& = \mathbf{A}+\Delta{A}
\end{aligned}
\end{equation}
with $\Delta{A}\in\mathbb{R}^{m\times m}$ denoting the error matrix introduced by the truncation and $\Delta_{ij}/m$ the $(i,j)$th element of $\Delta{A}$. It follows from Eq.(\ref{eq.13}) that $|\Delta_{ij}|\leq 2d\delta$ and $||\Delta{A}||_F\leq 2d\delta$, that is,
$||\Delta{A}||_F=\big|\big|A-\mathbf{A}\big|\big|_F \leq 2d\delta$.
In order to ensure that the error of preparing the interpolation matrix $\mathbf{A}$ is bounded by $\epsilon_{A}$, we select $\delta=\epsilon_{A}/(2d)$ such that
$$||\Delta{A}||_F\leq \epsilon_{A}.$$
Consequently, the total cost to prepare $A$ is $O\left(md\log(1/\delta)\right)=O\left(md\log \left(d/\epsilon_{A}\right)\right)$.

Note that $A$ is also a symmetric matrix and we assume that $\epsilon_{A}$ is small enough such that $A$ keeps positive definite. Without confusion, we also refer $A$ as the interpolation matrix. To summarize, we can quantum mechanically prepare the interpolation matrix with accuracy $\epsilon_{A}$ and computational complexity $O\left(md\log (d/\epsilon_{A})\right)$. Compared with the classical methods, a quadratic speedup is achieved in the number $m$ of the data points.

\subsubsection{$\mathbf{Interpolation \; matrix \; exponentiation}$}
Having derived the interpolation matrix, we want to find a quantum matrix inversion algorithm. For quantum mechanically computing a matrix inverse such as $A^{-1}$, one need to be able to enact $e^{-iAt}$ efficiently for some $t\geq0$. However, since $A$ defined by Eq.(\ref{eq.15}) is not sparse, techniques developed in \cite{berry2007efficient} cannot be used. For the exponentiation of nonsparse symmetric or Hermitian matrices, we use the strategy developed in \cite{lloyd2014quantum} to the present problem. The efficient preparation and exponentiation of the interpolation matrix, which appears in many scattered data interpolation problems, potentially enables a wide range of quantum scattered data-fitting algorithms.

For any density matrix $\rho$, applying the algorithm of \cite{lloyd2014quantum} we obtain
\begin{equation}\label{eq.16}
\begin{aligned}
{\rm tr}_1\{e^{-iS\Delta t}A\otimes\rho e^{iS\Delta t}\}
	& = \rho-i\Delta t[A,\rho]+O(\Delta t^2)\\
	& \approx e^{-iA\Delta t}\rho e^{iA\Delta t},
\end{aligned}
\end{equation}
where ${\rm tr}_1$ is the partial trace over the first variable and $S=\sum_{j,k=1}^m\ket{j}\bra{k}\otimes \ket{k}\bra{j}$ is the swap matrix of dimension $m^2\times m^2$. $S$ is a 1-sparse matrix and  $e^{-iS\Delta t}$ can be performed efficiently in negligible time $\widetilde{O}\left(\log\left(m\right) \Delta t\right)$ \cite{berry2007efficient}.
Eq.(\ref{eq.16}) shows that enacting $e^{-iA\Delta t}$ is possible with error $O(\Delta t^2)$. Repeated application of Eq.(\ref{eq.16}) with $l$ copies of $A$ allows one to construct $e^{-iAt}\rho e^{iAt}$, i.e.,
\begin{equation*}
	\begin{split}
		& {\rm tr}_1\Big{\{}e^{-iSt/l}A\otimes\cdots{\rm tr}_1\{e^{-iSt/l}A\otimes\rho e^{iSt/l}\}\cdots e^{iSt/l}\Big{\}}_l\\
		&= e^{-iAt}\rho e^{iAt}+O(t^2/l)
	\end{split}
\end{equation*}
with $t=l\Delta t$. To simulate $e^{-iAt}$ to accuracy $\epsilon_{E}$, it requires $l=O\left(t^2/\epsilon_E\right)$ steps. That is, $O\left(t^2/\epsilon_E\right)$ copies of $A$ allow us to construct the unitary operator $e^{-iAt}$ up to error $\epsilon_{E}$ \cite{lloyd2014quantum}. Taking into account the preparation of $A$ in time $O\left(md\log\left(d/\epsilon_{A}\right)\right)$, the runtime to simulate $e^{-iAt}$ is thus
$$
O\left(t^2\epsilon_{E}^{-1}md\log\left(d\epsilon_{A}^{-1}\right)\right)
$$
within error $\epsilon_{E}$.

The construction here shows that multiple copies of $A$ can be used to implement the unitary operator $e^{-iAt}$, where $A$ can be regard as an energy operator or Hamiltonian. It follows that the eigenvalues and eigenvectors of $A$ can be represented in quantum form and further used to perform a quantum matrix inversion.

\subsubsection{$\mathbf{Quantum \; globally \; supported \; RBF \; algorithm \; implementation}$}
By the interpolation matrix preparation and exponentiation techniques developed above, we present a quantum RBF algorithm in this subsection. The aim is to generate a quantum state $\ket{\bm c^{(exact)}}=\ket{\mathbf{A}^{-1}\bm y}$ corresponding to the solution of the linear system (\ref{eq.2}) with $\mathbf{A}$ defined by (\ref{eq.5}) and then compute the value at a new given data via $f$ in the form of Eq.(\ref{eq.4}). However,
since we can efficiently prepare the interpolation matrix $A$ defined by Eq.(\ref{eq.15}) and enact $e^{-iAt}$, we would like to solve the linear system
\begin{equation}\label{eq.17}
	A\bm c=\bm y
\end{equation}
and find a quantum state $\ket{\bm c}=\ket{A^{-1}\bm y}$ corresponding to the solution of this linear system. The right hand side $\bm y$ has been defined in (\ref{eq.2}). From Eq.(\ref{eq.15}), we know that $A$ is an approximation of $\mathbf{A}$, that is, $A=\mathbf{A}+\Delta A$
with $||\Delta A||_F\leq \epsilon_A$. Therefore, $\bm c$ is approximate to $\bm c^{(exact)}$ defined in (\ref{eq.2}), and its quantum state
$\ket{\bm c}$ is approximate to $\ket{\bm c^{(exact)}}$. We will analyze the error in the next subsection.

With the ability to construct the unitary transformation $e^{-iAt}$, the quantum matrix inversion technique developed in  \cite{harrow2009quantum} then allows us to obtain the quantum state $\ket{\bm c}$. The value at a new given data will be determined by the probability of a swap test. We now give a complete quantum globally supported RBF algorithm.

Denote by $\ket{u_j}$ the unit eigenvectors of $A$, and by $\lambda_j$ the corresponding eigenvalues. Besides, we use $\lambda_{\max}$ and $\lambda_{\min}$ to represent the largest and smallest eigenvalues of $A$. It follows from Lemma \ref{eq.2} that
\begin{equation}\label{eq.18}
|\lambda_{\max}-\Lambda_{\max}|\leq\epsilon_{A},\ |\lambda_{\min}-\Lambda_{\min}|\leq\epsilon_{A},
\end{equation}
where $\Lambda_{\max}$ and $\Lambda_{\min}$ are the largest and smallest eigenvalues of $\mathbf{A}$. The quantum state $\ket{\bm y}$ of the right hand side in (\ref{eq.17}) can be decomposed in the eigenvectors $\ket{u_j}$ of $A$, that is,
\begin{equation}\label{eq.19}
\ket{\bm y}=\sum_{j=1}^{m}\beta_j\ket{u_j}
\end{equation}
with $\beta_j=\braket{u_j|{\bm y}}$. For application of the quantum matrix inversion algorithm, we employ $e^{-iAt}$ conditionally in phase estimation. With a register initialized to $\ket{0}$ for storing an approximation of eigenvalues, phase estimation generates a state which is close to the ideal state storing the respective eigenvalues
$$
\sum\limits_{j=1}^{m}\beta_j\ket{u_j}\ket{\lambda_j}.
$$
The next step inverts the eigenvalue. Adding an ancilla qubit, performing a controlled rotation and uncomputing the eigenvalue register yields
\begin{equation}\label{eq.20}
\sum\limits_{j=1}^{m}\beta_j\ket{u_j}\left(\sqrt{1-\frac{C^2}{\lambda_j^2}}\ket{0}+\frac{C}{\lambda_j}\ket{1}\right),
\end{equation}
where $C=O(\lambda_{\min})$ is a constant such that the rotation can be defined and the probability of measuring the ancilla to be 1 is not greater than 1. To finish the inversion, measure the last qubit. Conditioned on seeing 1, we have the state $\ket{\bm c}$, that is,
$$
\sqrt{\frac{1}{\sum_{j=1}^{m}C^2|\beta_j|^2/|\lambda_j|^2}}\sum\limits_{j=1}^{m}\beta_j\frac{C}{\lambda_j}\ket{u_j}
$$
corresponding to $\sum_{j=1}^{m}\beta_j\lambda_j^{-1}\ket{u_j}$ up to a normalization. Note that the solution $\bm c$ in (\ref{eq.17}) can be expressed as
$$
\bm c=A^{-1}\bm y= A^{-1}\left(||\bm y||\cdot\ket{\bm y}\right)=||\bm y||\sum_{j=1}^{m}\beta_j\lambda_j^{-1}\ket{u_j}
$$
due to $\bm y = ||\bm y||\cdot\ket{\bm y}$ and Eq.(\ref{eq.19}), where $||\bm y||$ is given according to our assumption. We can determine the normalization factor
\begin{equation}\label{eq.21}
F\equiv\sqrt{\sum_{j=1}^{m}C^2|\beta_j|^2/|\lambda_j|^2}=\frac{C}{||\bm y||}||\bm c||
\end{equation}
from the probability of obtaining 1 in Eq.(\ref{eq.20}). We know from statistical sampling that $O(1/\epsilon_{F}^2)$ samples are required to estimate the normalization factor $F$ within error $\epsilon_{F}$. Consequently, we can obtain the quantum state  $\ket{\bm c}$ corresponding to the solution of the linear system (\ref{eq.17}), and have an estimate on the norm $||\bm c||$ of the solution from Eq.(\ref{eq.21}).

To evaluate on a new given data $\bm x$ via $f$ defined by Eq.(\ref{eq.4}), we apply a swap test to estimate the inner product between $\ket{\bm c}$ and $\ket{\Phi(\bm x)}$, where $\ket{\Phi(\bm x)}$ is the quantum state of vector
$$
\Phi(\bm x)=[e^{-||\bm x-\bm{x}^{(1)}||^2/2\sigma^2}, e^{-||\bm x-\bm{x}^{(2)}||^2/2\sigma^2},\dots, e^{-||\bm x-\bm{x}^{(m)}||^2/2\sigma^2}]^T.
$$
This requires the preparation of $\ket{\Phi(\bm x)}$. To this end we evaluate all the $m$ entries of $\Phi(\bm x)$ and $||\Phi(\bm x)||$ with $O(md)$ time by classical methods. Then produce the quantum state $\ket{\Phi(\bm x)}$ by a standard technique with QRAM. Using an ancilla, construct the state
$$
\frac{1}{\sqrt{2}}\left(\ket{0}\ket{\bm c}\ket{\Phi(\bm x)}+\ket{1}\ket{\Phi(\bm x)}\ket{\bm c}\right).
$$
Then measure the ancilla in the state $\ket{+}=1/\sqrt{2}(\ket{0}+\ket{1}).$ The measurement has the success probability
\begin{equation}\label{eq.22}
\begin{aligned}
p
&=\big|\big|(\bra{+}\otimes I)\frac{1}{\sqrt{2}}\left(\ket{0}\ket{\bm c}\ket{\Phi(\bm x)}+\ket{1}\ket{\Phi(\bm x)}\ket{\bm c}\right)\big|\big|^2\\
&=\frac{1}{2}+\frac{1}{2}\big|\braket{\bm c|\Phi(\bm x)}\big|^2,
\end{aligned}
\end{equation}
from which we can estimate $\braket{\bm c|\Phi(\bm x)}$. Similarly, using statistical sampling we can obtain $p$ to accuracy $\epsilon_p$ by iterating $O(1/\epsilon_p^2)$ times. Thus, the inner product $\braket{\bm c|\Phi(\bm x)}$ can be estimated with error $O(\epsilon_p)$. Given that $\ket{\bm c}$ and $\ket{\Phi(\bm x)}$ are quantum states, a multiplier is required to estimate the inner product between vectors $\bm c$ and $\Phi(\bm x)$. It is easy to find that this multiplier is of the form $||\bm c||\cdot||\Phi(\bm x)||$, where $||\bm c||$ and $||\Phi(\bm x)||$ can be evaluated. Therefore, by multiplying the multiplier to $\braket{\bm c|\Phi(\bm x)}$ we can obtain the inner product between the vectors $\bm c$ and $\Phi(\bm x)$, which is an approximation of the value at the new given data $\bm x$, i.e.,
$$
f(\bm{x})\approx ||\bm c||\cdot||\Phi(\bm x)||\braket{\bm c|\Phi(\bm x)}.
$$

Finally, we summarize the procedure for the quantum globally supported RBF method in Algorithm \ref{algorithm.1}. Note that $\ket{\bm c}$ is the ideal state we wish to obtain from the quantum matrix inversion algorithm and it corresponds to the solution of the linear system (\ref{eq.17}). In fact, due to the error which mainly comes from the quantum phase estimation, we obtain a quantum state approximate to $\ket{\bm c}$. We also obtain an approximation of the value at a new given data. The detailed error analysis is contained in the following subsection.
\begin{algorithm}[t]
\caption{Quantum globally supported RBF algorithm.}
\label{algorithm.1}
\begin{algorithmic}[l]
\STATE {1. Prepare the interpolation matrix $A$.}
\STATE {2. Construct the unitary operator $e^{-iAt}$ by repeated copies of $A$.}
\STATE {3. Perform the matrix inversion algorithm to generate a quantum state approximate to $\ket{\bm c}=\ket{A^{-1}\bm y}$.}
\STATE {4. Evaluate on a new given data $\bm x$ by the swap test between $\ket{\bm c}$ and $\ket{\Phi(\bm x)}$.}	
\end{algorithmic}	
\end{algorithm}

\textbf{Remark 2.} \ By statistically sampling the outcomes from many instances of the swap test, the value $p$ defined by Eq.(\ref{eq.22}) can be learned. $p$ can be approximated by the sample mean of this distribution. From estimate of the standard deviation of the mean, $O(1/\epsilon_p^2)$ samples are required to estimate the mean within error $\epsilon_p$. The value of the normalization factor defined by Eq.(\ref{eq.21}) can be estimated in a similar way.

\subsubsection{$\mathbf{Error \; analysis \; and  \; runtime \; estimate}$}
We continue with error analysis and runtime estimate for this quantum algorithm. There are two main error sources. One comes from the truncation error in Subsection 4.1.1, and the other one mainly comes from the quantum phase estimation in Subsection 4.1.3. For clarity, denote by $\ket{\bm c^{(exact)}}=\ket{\mathbf{A}^{-1}\bm y}$ the exact state we wish to obtain corresponding to the solution of the linear system (\ref{eq.2}) with the coefficient matrix in the form of (\ref{eq.5}), by $\ket{\bm c}=\ket{A^{-1}\bm y}$ the ideal state corresponding to the solution of the linear system (\ref{eq.17}) provided that everything is exact for the quantum matrix inversion algorithm, and by $\ket{\bm c^{(actual)}}$ the actual output state taking into account the error. The task for the error analysis is to find the distance between $\ket{\bm c^{(exact)}}$ and $\ket{\bm c^{(actual)}}$, i.e., $\big|\big|\ket{\bm c^{(exact)}}-\ket{\bm c^{(actual)}}\big|\big|$.

By the property of vector norm, it is easy to derive
$$
\big|\big|\ket{\bm c^{(exact)}} - \ket{\bm c^{(actual)}}\big|\big|
\leq \big|\big|\ket{\bm c^{(exact)}}-\ket{\bm{c}}\big|\big|+ \big|\big|\ket{\bm c}-\ket{\bm c^{(actual)}}\big|\big|,
$$
in which the first item of the right hand side defines the error introduced by the truncation and the second defines the error from the quantum matrix inversion algorithm. We analyze them one by one as follows.

We first discuss the error between $\ket{\bm c^{(exact)}}$ and $\ket{\bm c}$, which is mainly induced by an approximation of the interpolation matrix (\ref{eq.5}). Ideally, we have $\mathbf{A}$ from Eq.(\ref{eq.14}) to perform the quantum matrix inversion. In fact, we apply $A=\mathbf{A}+\Delta A$ defined by Eq.(\ref{eq.15}) with $||\Delta A||_F\leq \epsilon_{A}$ to design our algorithm. We show the influence on the result of this algorithm. Set $\gamma\equiv \Big|\Big|\mathbf{A}^{-1}\Delta A\Big|\Big|_2<1$. It follows from Lemma \ref{lemma.1} that
\begin{equation}
\begin{aligned}
\big|\big|\ket{\bm c^{(exact)}}-\ket{\bm{c}}\big|\big|
& = \Bigg|\Bigg|\frac{\bm c^{(exact)}}{||\bm c^{(exact)}||}-\frac{\bm{c}}{||\bm{c}||}\Bigg|\Bigg|\\
& = \Bigg|\Bigg|\frac{\mathbf{A}^{-1}\bm y}{||\mathbf{A}^{-1}\bm y||}-\frac{A^{-1}\bm y}{||A^{-1}\bm y||}\Bigg|\Bigg|\\
& \leq \Bigg|\Bigg|\frac{\mathbf{A}^{-1}\bm y}{||\mathbf{A}^{-1}\bm y||}-\frac{A^{-1}\bm y}{||\mathbf{A}^{-1}\bm y||}\Bigg|\Bigg|+\Bigg|\Bigg|\frac{A^{-1}\bm y}{||\mathbf{A}^{-1}\bm y||}-\frac{A^{-1}\bm y}{||A^{-1}\bm y||}\Bigg|\Bigg|\\
& \leq \frac{2||A^{-1}\bm y-\mathbf{A}^{-1}\bm y||}{||\mathbf{A}^{-1}\bm y||}\\
& \leq 2||A^{-1}-\mathbf{A}^{-1}||_2\cdot ||\mathbf{A}||_2\\
& \leq \frac{2||\Delta{A}||_2\cdot ||\mathbf{A}^{-1}||^2_2\cdot ||\mathbf{A}||_2}{1-\gamma}\\
& \leq \frac{2\epsilon_{A}}{(1-\gamma)\Lambda_{\max}}\cdot \kappa^2,
\end{aligned}
\end{equation}
in which the last inequality uses the fact that $||\Delta{A}||_2\leq||\Delta{A}||_F\leq\epsilon_{A}$ and  $\kappa=||\mathbf{A}||_2\cdot||\mathbf{A}^{-1}||_2=\Lambda_{\max}/\Lambda_{\min}$. Here, $\Lambda_{\max}$ and $\Lambda_{\min}$ represent the largest and smallest eigenvalues of $\mathbf{A}$ respectively. In order to make $||\ket{\bm c^{(exact)}}-\ket{\bm{c}}||\leq\epsilon_{\bm c}/2$, let
$$
\frac{2\epsilon_{A}}{(1-\gamma)\Lambda_{\max}}\cdot \kappa^2\leq\frac{\epsilon_{\bm c}}{2}.
$$
Then we have
$$
\epsilon_{A}\leq\frac{1}{4}(1-\gamma)\Lambda_{\max}\cdot\epsilon_{\bm c}\kappa^{-2}
\leq\epsilon_{\bm c}\kappa^{-2}
$$
by the fact that $\Lambda_{\max}\leq 1$ from Eq.(\ref{eq.6}). Hence, we can choose $\epsilon_{A}=O(\epsilon_{\bm c}/\kappa^2)$ such that
$$
\big|\big|\ket{\bm c^{(exact)}}-\ket{\bm{c}}\big|\big|\leq\frac{\epsilon_{\bm c}}{2}.
$$

We then analyze the error from the quantum matrix inversion algorithm in Subsection 4.1.3, in which the dominate source of error is phase estimation. Let $t_0$ be the total evolution time determining the error of phase estimation \cite{harrow2009quantum}. The error for estimating $\lambda$ is $O(1/t_0)$ in this step, with $\lambda$ representing the eigenvalues of $A$. The relative error of $\lambda^{-1}$ is given by $O(1/(\lambda t_0))$. Since $\lambda\geq \lambda_{\min}$, we can take $t_0=O\left(1/\left(\lambda_{\min}\epsilon_{\bm c}\right)\right)$ such that
$$
\big|\big|\ket{\bm{c}}-\ket{\bm{c}^{(actual)}}\big|\big|\leq \frac{\epsilon_{\bm c}}{2}.
$$

From Subsection 4.1.2, we know that it requires $O(t_0^2\epsilon_{E}^{-1}md\log (d\epsilon_{A}^{-1}))$ time to perform the quantum matrix inversion. With $\epsilon_{A}=O(\epsilon_{\bm c}/\kappa^2)$ and $t_0=O\left(1/\left(\lambda_{\min}\epsilon_{\bm c}\right)\right)$, the runtime is thus
$O\left(\lambda_{\rm min}^{-2}\epsilon_{\bm c}^{-2}\epsilon_{E}^{-1}md\log(d\kappa^{2}\epsilon_{\bm c}^{-1})\right)$
to obtain a quantum state corresponding to the solution of the linear system (\ref{eq.17}) up to accuracy $\epsilon_{\bm c}$, that is,
$$
\big|\big|\ket{\bm c^{(exact)}} - \ket{\bm{c}^{(actual)}}\big|\big|\leq\epsilon_{\bm c}.
$$

Finally, we consider the success probability of the post-selection process on Eq.(\ref{eq.20}). Since $C=O(\lambda_{\min})$ and $\lambda\leq\lambda_{\max}\leq 1+\epsilon_{A}$ from Eq.(\ref{eq.6}) and (\ref{eq.18}), this probability is at least $\Omega(\lambda_{\min}^2)$. Using amplitude amplification \cite{brassard2002quantum}, we find that $O(1/\lambda_{\rm min})$ repetitions are sufficient. Furthermore, we analyze the upper bound of $1/\lambda_{\min}$. It follows from Eq.(\ref{eq.18}) that
$$
\frac{1}{\lambda_{\min}}\leq \frac{1}{\Lambda_{\min}-\epsilon_{A}}\leq\frac{2}{\Lambda_{\min}}=\kappa\cdot \frac{2}{\Lambda_{\max}},
$$
that is, $1/\lambda_{\rm min}=O(\kappa)$. Putting this all together, choosing $\epsilon=\min\{\epsilon_E, \epsilon_{\bm c}\}$, we obtain the total runtime of
$$
O(\kappa^3\epsilon^{-3}md\log(d\kappa^2\epsilon^{-1}))
$$
to derive the final state $\ket{\bm c^{(exact)}}$ with precision $\epsilon$.

The error for computing the value at a new given data comes from the estimation on the norm $||\bm c||$ from Eq.(\ref{eq.21}) and the inner product $\braket{\bm c|\Phi(\bm x)}$ from
Eq.(\ref{eq.22}) on the one hand and an approximation of the quantum state $\ket{\bm c^{(exact)}}$ on the other hand. Choosing $\epsilon_{F}=O(\epsilon)$ and $\epsilon_p=O(\epsilon)$ with $\epsilon_F$ and $\epsilon_p$ representing the error for evaluating $F$ and $p$ defined by Eq.(\ref{eq.21}) and (\ref{eq.22}) respectively, it is not hard to derive the value at the new given data with error $O(\epsilon)$. Taking into account the cost $O(1/\epsilon^{2})$ for obtaining $||\bm c||$ and $\braket{\bm c|\Phi(\bm x)}$ with error $O(\epsilon)$, the total runtime for evaluating on a new given data is
$$
O(\kappa^3\epsilon^{-5}md\log(d\kappa^2\epsilon^{-1}))
$$
with error $O(\epsilon)$.

To summarize, we find a quantum RBF method for global scattered data interpolation. Table \ref{table.2} compares the performance of our approach with the best classical method in terms of the generator and solver, where the generator is used for generating the interpolation matrix and the solver is used for solving the linear system. It shows that our quantum algorithm improves the dependence on the number $m$ of data, which scales linearly in $m$, achieving a quadratic speedup over the best classical method. Specially, our algorithm could be a subroutine in a larger quantum algorithm where the interpolation matrix can be produced by some other methods. We summarize the results into a theorem as follows.
\begin{table}[h]
	\centering
	\caption{\label{table.2} Comparison about the time complexity of classical and quantum algorithms. Generator: generating the interpolation matrix. Solver: solving the linear system. $T\equiv md\log(d/\epsilon_A)$.}
	\begin{tabular}{lll}
		\hline\hline
		Algorithm & Generator & Solver\\
		\hline
		Classical & $O(m^2d)$ & $O\left(m^2\sqrt{\kappa}\log(1/\epsilon)\right)$ \cite{shewchuk1994introduction}\\
		Quantum & $O(T)$  & $O(\kappa^3\epsilon^{-3}T)$\\
		\hline\hline		
	\end{tabular}
\end{table}
\newtheorem{theorem}{Theorem}
\begin{theorem}
For Problem \ref{problem.1}, we find a quantum algorithm to interpolate the scattered data using the Gaussian RBF. The cost for obtaining a quantum state corresponding to the coefficient vector is
$$
O(\kappa^3\epsilon^{-3}md\log\left(d\kappa^2\epsilon^{-1})\right)
$$
up to accuracy $\epsilon$, and the cost for calculating the value at a new given data by the swap test is
$$
O(\kappa^3\epsilon^{-5}md\log(d\kappa^2\epsilon^{-1}))
$$
with error $O(\epsilon)$, where $\kappa$ is the condition number of the interpolation matrix $\mathbf{A}$ defined by (\ref{eq.5}).	
\end{theorem}

\subsection{$\mathbf{Quantum \; compactly \; supported \; RBF \; method}$}
So far, we have designed a quantum algorithm for the globally supported RBF interpolation problem. This approach has the disadvantage of being global in the sense that every data point has some influence on the solution function at any point. It would be more interesting to allow only the nearest neighbors of the evaluation point to influence the approximate value. In numerical analysis, the concept of locally supported basis functions is of great importance. The general advantage of compactly supported basis functions lies on a sparse interpolation matrix and the possibility of a fast evaluation of the interpolant. In this subsection, we concentrate on positive definite radial functions with compact support for the scattered data interpolation problem.

There are many ways to construct positive definite radial functions with compact support \cite{wendland1995piecewise, wu1995compactly, gneiting2002compactly}. Probably the most popular family of compactly supported radial functions was constructed by Wendland in \cite{wendland1995piecewise}. Wendland has solved for the minimum-degree polynomial solution for compactly supported RBFs which guarantee positive-definiteness of the matrix. All of the radial functions of Wendland have the form
\begin{equation}\label{eq.24}
\phi(r)=\left\{
	\begin{array}{lcl}
		(1-r)_+^lq(r), &  & 0\leq r \leq 1,\\
		0, & & r>1,
	\end{array}
\right.
\end{equation}
with a polynomial $q$. Here the cutoff function $(\cdot)_+$ is defined by
$$
w_+=\left\{
\begin{array}{lcl}
	w, &  & w\geq 0,\\
	0, & & w<0.
\end{array}
\right.
$$
We also point out that $(1-r)_+^l$ in Eq.(\ref{eq.24}) is interpreted as $((1-r)_+)^l$, that is, we first apply the cutoff function, and then the power. For various degrees of space dimension ($d$) and continuity ($C^k$) of the interpolated function, Wendland \cite{wendland1995piecewise} has derived the following compactly supported RBFs in Table \ref{table.3}.
\begin{table}[h]
	\centering
	\caption{\label{table.3}Compactly supported functions of minimal degree.}
	\begin{tabular}{llc}\hline\hline
		Space dimension & Function & Smoothness\\ \hline
		$d=1$ & $(1-r)_+$ & $C^0$\\
		 & $(1-r)_+^3(3r+1)$ & $C^2$\\
		 & $(1-r)_+^5(8r^2+5r+1)$ & $C^4$\\
		 \hline
		$d=3$ & $(1-r)_+^2$ & $C^0$\\
		 & $(1-r)_+^4(4r+1)$ & $C^2$\\
		 & $(1-r)_+^6(35r^2+18r+3)$ & $C^4$\\
		 & $(1-r)_+^8(32r^3+25r^2+8r+1)$ & $C^6$\\
		 \hline
		$d=5$ & $(1-r)_+^3$ & $C^0$\\
		& $(1-r)_+^5(5r+1)$ & $C^2$\\
		& $(1-r)_+^7(16r^2+7r+1)$ & $C^4$\\
		\hline\hline		
	\end{tabular}
\end{table}
These functions have radius of support equal to 1. Scaling of the functions, that is, $\phi(r/\alpha)$, allows any desired radius of support $\alpha$. The interpolation matrix corresponding to this family of RBFs will be sparse by scaling the support $\alpha$ of the basis function $\phi_\alpha(r)\equiv\phi(r/\alpha)$ appropriately.

In the following, we will present a quantum compactly supported RBF algorithm by means of this family of compactly supported RBFs for the scattered data interpolation problem. In this case, the approximation function $f$ reads
\begin{equation}\label{eq.25}
	f(\bm{x}) = \sum\limits_{j=1}^{m}c_j \phi(||\bm x-\bm x^{(j)}||/\alpha),
\end{equation}
where $\alpha$ is a scaling parameter. The coefficients $c_j$ are found by solving the linear system (\ref{eq.1}) with the interpolation matrix in the form
\begin{equation}\label{eq.26}
	A =
	\left[
	\begin{array}{cccc}
		\phi(||\bm{x}^{(1)} - \bm{x}^{(1)}||/\alpha) & \phi(||\bm{x}^{(1)} - \bm{x}^{(2)}||/\alpha) & \cdots & \phi(||\bm{x}^{(1)} - \bm{x}^{(m)}||/\alpha)\\
		\phi(||\bm{x}^{(2)} - \bm{x}^{(1)}||/\alpha) & \phi(||\bm{x}^{(2)} - \bm{x}^{(2)}||/\alpha) & \cdots & \phi(||\bm{x}^{(2)} - \bm{x}^{(m)}||/\alpha)\\
		\vdots & \vdots & \ddots & \vdots \\
		\phi(||\bm{x}^{(m)} - \bm{x}^{(1)}||/\alpha) & \phi(||\bm{x}^{(m)} - \bm{x}^{(2)}||/\alpha) & \cdots & \phi(||\bm{x}^{(m)} - \bm{x}^{(m)}||/\alpha)
	\end{array}
	\right].
\end{equation}

\subsubsection{$\mathbf{Quantum \; compactly \; supported \; RBF \; algorithm \; implementation}$}
Since the interpolation matrix $A$ defined by (\ref{eq.26}) is sparse, we can use the HHL algorithm to solve the linear system. We assume that there are no more than $s$ nonzero elements in any row or column of $A$. The elements of $A$ defined by Eq.(\ref{eq.24}) are not directly available to us. We require an oracle to compute the nonzero entries
\begin{equation}\label{eq.27}
\mathcal{P}_A: \ket{i, j}\ket{0} \mapsto \ket{i,j}\ket{A_{ij}}
\end{equation}
for any $i,j\in\{1,2,\dots,m\}$, where $A_{ij}$ is the binary representation of the $(i,j)$th element of $A$ in (\ref{eq.26}).

We first bulid this oracle and show how it is possible to compute the values of the nonzero entries of $A$. Using an ancilla, construct the state by calling the data sites oracle
\begin{equation}
\ket{\phi_{ij}}=\frac{1}{\sqrt{||\bm x^{(i)}||^2+||\bm x^{(j)}||^2}}\left(\big|\big|\bm x^{(i)}\big|\big|\ket{+}\ket{\bm x^{(i)}}-\big|\big|\bm x^{(j)}\big|\big|\ket{-}\ket{\bm x^{(j)}}\right)
\end{equation}
for any $i,j\in\{1,2,\dots,m\}$ in quantum parallel, which takes $O(\log d)$ time. Here, $\ket{+}=1/\sqrt{2}(\ket{0}+\ket{1})$ and $\ket{-}=1/\sqrt{2}(\ket{0}-\ket{1})$ can be easily produced by the Hadamard transform. Starting from the initial state as $\ket{i,j}\ket{0}$, we can then obtain
$$
\ket{i,j}\ket{0} \mapsto \ket{i,j}\ket{\phi_{ij}}.
$$
Note that the amplitude of $\ket{0}$ of the first register of  $\ket{\phi_{ij}}$ is
$$
||\bm x^{(i)}-\bm x^{(j)}||\Big/{\sqrt{2\left(||\bm x^{(i)}||^2+||\bm x^{(j)}||^2\right)}}
$$
corresponding to the distance between the vectors $\bm x^{(i)}$ and $\bm x^{(j)}$. Since $||\bm x^{(i)}||$ for $i\in\{1,2,\dots,m\}$ are known by the assumption, an approximation of $||\bm x^{(i)}-\bm x^{(j)}||$ can be obtained with the help of the amplitude estimation \cite{brassard2002quantum}. With this approximation storing in an ancilla register, we obtain a state which is close to the ideal state storing the distance between $\bm x^{(i)}$ and $\bm x^{(j)}$ in the ancilla
\begin{equation}\label{eq.29}
\ket{i,j}\ket{\phi_{ij}}\ket{||\bm x^{(i)}-\bm x^{(j)}||}.
\end{equation}
Let $\mathcal{O}_{\phi}$ be an oracle to compute the value of $\phi$ defined by Eq.(\ref{eq.24}). Adding an ancilla qubit $\ket{0}$ and applying the oracle $\mathcal{O}_\phi$, we get
\begin{equation}\label{eq.30}
\ket{i,j}\ket{\phi_{ij}}\ket{||\bm x^{(i)}-\bm x^{(j)}||}\ket{\phi(||\bm x^{(i)}-\bm x^{(j)}||/\alpha)}.
\end{equation}
Uncomputing the registers of $\ket{\phi_{ij}}$ and $\ket{||\bm x^{(i)}-\bm x^{(j)}||}$ in Eq.(\ref{eq.30}) we obtain
$$
\ket{i,j}\ket{\phi(||\bm x^{(i)}-\bm x^{(j)}||/\alpha)},
$$
that is,
$$
\ket{i,j}\ket{A_{ij}}.
$$
Denote the above steps by $\mathcal{P}_A$. We have
\begin{equation}
\mathcal{P}_A\ket{i,j}\ket{0}=\ket{i,j}\ket{A_{ij}}.
\end{equation}
Thus, the oracle $\mathcal{P}_A$ defined by Eq.(\ref{eq.27}) can be derived.

Besides, to perform the Hamiltonian simulation of the sparse matrix $A$, we also need an oracle to compute the locations of the nonzero entries. We assume that $\mathcal{P}_{\upsilon}$ allows us to perform the map
\begin{equation}
	\mathcal{P}_{\upsilon}: \ket{j, \mathcal{\ell}} \mapsto \ket{j,\upsilon(j,\mathcal{\ell})}
\end{equation}
for $j\in\{1,2,\dots,m\}$ and $\mathcal{\ell}\in\{1,2,\dots,s\}$, where $s$ represents the maximum number of nonzero entries in any row or column of $A$, and $\upsilon(j,\ell)$ is a function that returns the row index of the $\ell$th nonzero entry of the $j$th column. The assumption is reasonable since we can efficiently compute and locate the nonzero entries of $A$ by the oracle $\mathcal{P}_A$ we have built.

Having derived the oracle $\mathcal{P}_A$ and $\mathcal{P}_\upsilon$, we can efficiently simulate $e^{-iAt}$ for some $t\geq0$. We use the HHL
algorithm to obtain a quantum state $\ket{\bm c^{(compact)}}$ corresponding to the solution of the linear system $(\ref{eq.1})$ with the coefficient matrix in the form of (\ref{eq.26}). Moreover, the norm $||\bm c^{(compact)}||$ of the solution can be determined in a similar way like $||\bm c||$ in Eq.(\ref{eq.21}). The detailed procedure of the quantum matrix inversion and estimate on the norm of the solution can be referred to in Subsection 4.1.3. We omit it here.

For calculating the value at a new given data $\bm x$ via Eq.(\ref{eq.25}), we also apply a swap test between $\ket{\bm c^{(compact)}}$ and $\ket{\Phi(\bm x)}$, where $\ket{\Phi(\bm x)}$ is the quantum state for the vector
$$
\Phi(\bm x)=[\phi{(||\bm x-\bm{x}^{(1)}||/\alpha)}, \phi{(||\bm x-\bm{x}^{(2)}||/\alpha)},\dots, \phi{(||\bm x-\bm{x}^{(m)}||/\alpha)}]^T.
$$
We prepare $\ket{\Phi(\bm x)}$ in a similar way like the oracle $\mathcal{P}_A$. This gives
$$
\frac{1}{\sqrt{m}}\sum\limits_{j=1}^{m}\ket{j}\ket{\phi{(||\bm x-\bm{x}^{(j)}||/\alpha)}}
$$
with $m$ the number of data. Adding an ancilla qubit and applying control rotation on $\ket{\phi{(||\bm x-\bm{x}^{(j)}||/\alpha)}}$ yields
\begin{equation}{\label{eq.33}}
\frac{1}{\sqrt{m}}\sum\limits_{j=1}^{m}\ket{j}\ket{\phi{(||\bm x-\bm{x}^{(j)}||/\alpha)}}\left(\hat{C}\phi{(||\bm x-\bm{x}^{(j)}||/\alpha)}\ket{0}+\sqrt{1-\hat{C}^2\phi^2{(||\bm x-\bm{x}^{(j)}||/\alpha)}}\ket{1}\right)
\end{equation}
where $\hat{C}$ is a scaling factor such that the rotation can be defined and the probability of obtaining $0$  on the ancilla register is not greater than $1$. Uncomputing the register $\ket{\phi{(||\bm x-\bm{x}^{(j)}||/\alpha)}}$ and measuring the last qubit in the state $\ket{0}$, we obtain $\ket{\Phi(\bm x)}$ proportional to
$$
\sum\limits_{j=1}^{m}\phi{(||\bm x-\bm{x}^{(j)}||/\alpha)}\ket{j}.
$$
In the same way as estimating $||\bm c||$ from Eq.(\ref{eq.21}), we can estimate $||\Phi(\bm x)||$ from the probability of obtaining 0 on the ancilla register in Eq.(\ref{eq.33}). Furthermore, by the swap test between $\ket{\bm c^{(compact)}}$ and $\ket{\Phi(\bm x)}$, we can evaluate on the new given data $\bm x$. In this regard, we have already described it in Subsection 4.1.3 and will not go into details.

\subsubsection{$\mathbf{Error \; analysis \; and \; runtime \;  estimate}$}
We further analyze the error and discuss the runtime about this quantum algorithm. The computational complexity for implementing the oracle $\mathcal{P}_A$ mainly comes from amplitude estimation. The preparation of $\ket{\phi_{ij}}$ is efficient in time $O(\log d)$. By the quantum phase estimation algorithm \cite{kitaev1995quantum}, it takes $O(\epsilon^{-1}\log^2d)$ time for the amplitude estimation up to precision $\epsilon$. Taking into account the runtime $\widetilde{O}(\kappa^2s^2\epsilon^{-1}\log m)$ of the HHL algorithm, the total runtime is
$$
\widetilde{O}(\kappa^2s^2\epsilon^{-2}\log m\log^2d)
$$
to obtain a quantum state corresponding to the solution of the linear system up to accuracy $\epsilon$, where $\kappa$ is the condition number of $A$ in (\ref{eq.26}).

Considering the cost $O(1/\epsilon^2)$ for estimating the norms $||\bm c ^{(compact)}||$, $||\Phi(\bm x)||$ and the inner product between $\ket{\bm c^{(compact)}}$ and $\ket{\Phi(\bm x)}$ with error $O(\epsilon)$, the runtime for evaluating on a new given data via Eq.(\ref{eq.25}) is thus
$$
\widetilde{O}(\kappa^2s^2\epsilon^{-4}\log m\log^2d)
$$
with error $O(\epsilon)$.

In summary, we present a quantum compactly supported RBF algorithm for the scattered data interpolation problem, achieving an exponential speedup in the number (i.e., $m$) of data and the dimension (i.e., $d$) of input vectors over the classical methods. We summarize the results into the following theorem.
\begin{theorem}
	For Problem \ref{problem.1}, we find a quantum compactly supported RBF algorithm to efficiently interpolate the scattered data. The cost for obtaining a quantum state corresponding to the coefficient vector is
	$$
	\widetilde{O}(\kappa^2s^2\epsilon^{-2}\log m\log^2d)
	$$
	up to accuracy $\epsilon$, and the cost for calculating the value at a new given data by the swap test is
	$$
	\widetilde{O}(\kappa^2s^2\epsilon^{-4}\log m\log^2d)
	$$
	with error $O(\epsilon)$, where $\kappa$ is the condition number of the interpolation matrix $A$ defined by (\ref{eq.26}).	
\end{theorem}

\section{Conclusion}
In this paper, we have shown that the RBF method, an important method in the scattered data fitting, can be implemented quantum mechanically with improved algorithmic complexity over the classical counterpart. Using globally and locally supported RBFs, we present two quantum algorithms to deal with the problem. The idea of the coherent states and the matrix exponentiation technique allow the construction of the interpolation matrix and the quantum matrix inversion for the globally supported RBF method, achieving a quadratic speedup in the number of data. As for the compactly supported RBF method, our quantum algorithm achieves an exponential speedup in the dimension of input vectors and the number of data by building upon the HHL algorithm. Our results can also be extended to cases where the interpolation matrix is ill-conditioned, arising from the possibility of having a data site that has almost zero overlap with another data site. In this case, we will define a constant $\delta_{\rm{eff}}$ such that only the eigenvalues greater than $\delta_{\rm{eff}}$ are taken into account, essentially defining an effective condition number  $k_{\rm{eff}}$\cite{li2010effective}. Then the filtering procedure described in \cite{harrow2009quantum} is employed in the phase estimation using the $\delta_{\rm{eff}}$. Future work may concern more general cases such as finding a least squares solution if the values of the data contain noise, finding an optimal parameter provided that the scaling  parameter is unknown, etc.

\bibliographystyle{elsarticle-num}
\bibliography{QRBFReference}

\end{document}